\documentstyle[amssymb,12pt]{amsart}

\newtheorem{theorem}{Theorem}
\newtheorem{corollary}[theorem]{Corollary}
\newtheorem{definition}[theorem]{Definition}

\newcommand{\res}{\upharpoonright}

\begin{document}

\baselineskip=18pt

  \begin{center}
     {\Large\bf Compactness of Loeb Spaces}
  \end{center}

  \begin{center}
 Renling Jin\footnote{The research of the first author
was supported by NSF postdoctoral fellowship \#DMS-9508887. 
This research was started when the first author
spent a wonderful year as a visiting assistant professor
in University of Illinois-Urbana Champaign during 94-95. 
He is grateful to the logicians there.} 
\& Saharon Shelah\footnote{The research
of the second author was supported by The Israel Science Foundation
administered by The Israel Academy of Sciences and Humanities. 
This paper is number 613 on the second author's publication list.}

  \end{center}

  \bigskip

  \begin{quote}

    \centerline{Abstract}

    \small  
In this paper we show that the compactness of a Loeb space depends
on its cardinality, the nonstandard universe 
it belongs to and the underlying model of set theory we live in.
In \S 1 we prove that Loeb spaces are compact under various assumptions,
and in \S 2 we prove that Loeb spaces are not compact under various other
assumptions. The results in \S 1 and \S 2 give a quite complete
answer to a question of D. Ross in [R1], [R2] and [R3].
 
\end{quote}

\section{Introduction}

In [R1] and [R2] D. Ross asked: {\em Are 
(bounded) Loeb measure spaces compact?}
J. Aldaz then, in [A], constructed a counterexample. But Aldaz's example
is atomic, while most of Loeb measure spaces people are interested 
are atomless. So Ross re-asked his question in [R3]: {\em Are atomless
Loeb measure spaces compact?} In this paper we answer the question.
Let's assume that all measure spaces mentioned throughout this paper 
are \underline{atomless} probability spaces.

Given a probability space $(\Omega,\Sigma,P)$. A subfamily
${\cal C}\subseteq\Sigma$ is called compact if for any ${\cal D}
\subseteq {\cal C}$, $\cal D$ has f.i.p. {\em i.e.} 
finite intersection property,
implies $\bigcap {\cal D}\not=\emptyset$.
We call a compact family $\cal C$ inner-regular on $\Omega$ if
for any $A\in\Sigma$

\[P(A)=\sup \{P(C):C\subseteq A\wedge C\in {\cal C}\}.\]

Ross called a probability space $(\Omega,\Sigma,P)$ compact if
there is an inner-regular compact family on it. Clearly,
Ross's definition of compactness is a generalization of Radon spaces
with no topology involved.
In fact, Ross proved in [R2] that a compact probability space is
essentially Radon, {\em i.e.} one can topologize the space so that
the resulting space is Radon.

Loeb measure spaces are important tools in nonstandard analysis
(see, for example, [AFHL] and [SB]). Ross proved in [R2] 
that every compact probability
space is the image, under a measure preserving transformation, of a Loeb
measure space. This shows, by a word of Ross,
some evidence that Loeb spaces
themselves may be compact.

In this paper we show that the compactness of a Loeb space
depends on its cardinality, the nonstandard universe it belongs to,
and even the underlying world of set theory we live in (suppose we
live in a transitive model of ZFC).

Throughout this paper we always denote $\cal M$ for our underlying
transitive model of set theory ZFC. We sometimes use $\cal N$ for
another transitive model of ZFC. If we make a statement without
mentioning a particular model, this statement is always assumed to be
relative to $\cal M$. Let $\Bbb N$ be the set of all
standard natural numbers. Using $\Bbb N$ as a set of urelements,
we construct the standard universe $(V,\in)$ by letting
\[V_0={\Bbb N},\;V_{n+1}=V_n\cup {\cal P}(V_n)\mbox{ and }
V=\bigcup_{n\in\omega}V_n.\]
A nonstandard universe $(^*V,\;^*\!\in)$ is the truncation,
at ${}^*\!\in$-rank $\omega$, of
an elementary extension of the standard universe such that
$^*{\Bbb N}\smallsetminus {\Bbb N}\not=\emptyset$.
We always assume the nonstandard universe ${}^*V$ 
we work within is $\omega_1$-saturated. In fact, $\omega_1$-saturation
is needed in Loeb measure construction. For any
set $S$ we use $|S|$ for its set theoretic cardinality. If $S$ is
an internal set (in $^*V$), then $^*|S|$ means the internal
cardinality of $S$. For any object $S$ in the 
standard universe we always denote
$^*S$ for its nonstandard version in $^*V$.
For example, if $\Omega$ is an internal set, then $^*{\cal P}(\Omega)$ 
denote the set of all internal subsets of $\Omega$. 
Let $\Sigma_0\subseteq \;{}^*{\cal P}(\Omega)$ be an internal algebra
and let $P:\Sigma_0\mapsto\;^*[0,1]$ be an internal finitely additive
probability measure. We call $(\Omega,\Sigma_0,P)$ an internal
probability space. Let $st:\;{}^*[0,1]\mapsto [0,1]$ be 
the standard part map. Then $(\Omega,\Sigma_0,st\!\circ\! P)$ is a
standard finitely additive probability space.
Then one can use $\Sigma_0$ to generate uniquely an 
$st\!\circ\! P$-complete $\sigma$-algebra $\Sigma$ and extend
$st\!\circ\! P$ uniquely to a standard complete countably additive
probability measure $L_P$. The space $(\Omega,\Sigma,L_P)$
is called a Loeb space generated by $(\Omega,\Sigma_0,P)$.
Let $H$ be a hyperfinite integer, 
{\em i.e.} $H\in\;{}^*{\Bbb N}\smallsetminus {\Bbb N}$. Let
$\Omega=\{1,2,\ldots,H\}$, let $\Sigma_0=\;{}^*{\cal P}(\Omega)$ and
let $P(A)=\;{}^*|A|/H$ for each $A\in\Sigma_0$. We call
$(\Omega,\Sigma_0,P)$ a hyperfinite internal space. The space
$(\Omega,\Sigma,L_P)$ generated by $(\Omega,\Sigma_0,P)$ as
above is called a hyperfinite Loeb
space. Hyperfinite Loeb spaces are
most useful among other Loeb spaces in nonstandard analysis.
For notational simplicity we prove the results only
for hyperfinite Loeb spaces in this paper.
From now on we denote the symbol $(\Omega,\Sigma,L_P)$ or just $\Omega$ 
without confusion,
exclusively for a \underline{hyperfinite Loeb space}.
Most of the results for hyperfinite Loeb spaces in this paper can be
easily generalized to Loeb spaces in the general sense
(see Fact 3 at the beginning of \S 1 and the comments after that).

In \S 1 we show when a hyperfinite Loeb space is compact. We prove the 
following results.

{\em Corollary 3}:  Suppose CH (Continuum Hypothesis) holds.
Suppose $|^*{\Bbb N}|=\omega_1$. Then every hyperfinite Loeb space in
$^*V$ is compact.

{\em Corollary 4}:  Suppose MA (Martin's Axiom) holds.
Suppose $^*V$ is $2^{\omega}$-saturated and 
$|^*{\Bbb N}|=2^{\omega}$. Then every hyperfinite Loeb space in $^*V$
is compact.

{\em Corollary 5}:  Suppose $\cal M$ is obtained by adding
$\kappa$ Cohen reals to a ZFC model $\cal N$ for some
$\kappa\geqslant (2^{\omega})^{\cal N}$
with $\kappa^{\omega}=\kappa$ in $\cal N$.
Suppose $|^*{\Bbb N}|=2^{\omega}$ (in $\cal M$ now). Then every
hyperfinite Loeb space in $^*V$ is compact.

{\em Corollary 7}:  Suppose $\cal M$ is same as in Corollary 5.
In $\cal M$ suppose $\lambda$ is a strong limit cardinal 
with $cf(\lambda)\leqslant\kappa$.
Suppose $|^*{\Bbb N}|=\lambda$ and $^*V$ satisfies the 
$\aleph_0$-special model axiom (see [R5] or [J] for the definition). 
Then every hyperfinite Loeb space in $^*V$ is compact.

{\em Theorem 8}:  Suppose $\kappa$ is a strong limit cardinal
with $cf(\kappa)=\omega$ and $\lambda=\kappa^+=2^{\kappa}$.
Suppose $^*V$ is $\lambda$-saturated and has cardinality $\lambda$.
Then every hyperfinite Loeb space in $^*V$ is compact.

In \S 2 we show when a hyperfinite Loeb space is not compact. We prove the
following results.

{\em Theorem 9}:  Suppose $\lambda$ is a regular cardinal such that
$\kappa^{\omega}<\lambda$ for any $\kappa<\lambda$. Suppose
$|\Omega|=\lambda$. Then $\Omega$ is not compact.

{\em Theorem 10}:  Suppose $\lambda$ is a strong limit cardinal,
$\kappa=cf(\lambda)$ and $\mu^{\omega}<\kappa$ for any $\mu<\kappa$. 
Suppose $|\Omega|=\lambda$. Then $\Omega$ is not compact.

{\em Theorem 11}:  Suppose $\cal M$ is obtained by adding 
$\kappa$ random reals to a ZFC model $\cal N$ for some
regular $\kappa>(2^{\omega})^{\cal N}$ with $\kappa^{\omega}=\kappa$. 
Suppose $|\Omega|=2^{\omega}$ (note $2^{\omega}=\kappa$ in $\cal M$). 
Then $\Omega$ is not compact.

{\em Theorem 12}:  Suppose $\cal M$ is obtained by adding
$\kappa$ random reals to a ZFC model $\cal N$ for some
regular $\kappa>\omega$. Suppose $\lambda$ is a strong limit cardinal 
such that $cf(\lambda)\leqslant\kappa$. Suppose 
$|\Omega|=\lambda$ (hence $cf(\lambda)>\omega$). 
Then $\Omega$ is not compact.

{\em Theorem 13}:  Let $\lambda>|V|$ and $\lambda^{\omega}
=\lambda$. Then there exists a $^*V$ such that $|\Omega|=\lambda$
and $\Omega$ is not compact for every $\Omega$ in $^*V$. 

\medskip

In this paper we write $\lambda,\kappa,\mu,\ldots$ for cardinals,
$\alpha,\beta,\gamma,\ldots$ for ordinals and
$k,m,n,\ldots$ for natural numbers. We write
$\lambda^{\kappa}\,(\lambda^{<\kappa})$ for cardinal exponents
and $^{\kappa}\lambda\,(^{<\kappa}\lambda)$ for sets of functions.
For any set $S$ we write $[S]^{\lambda}$ for the set of all
subsets of $S$ with cardinality $\lambda$.
For any set $S$ we write $(^S2,\Sigma(^S2),\nu_S)$ for the 
complete probability space generated by all Baire sets of $^S2$
such that for any finite $S_0\subseteq S$ and any 
$\tau\in\; ^{S_0}2$, $\nu_S([\tau])=2^{-|S_0|}$, where $[\tau]=
\{f\in^S\!2:f\!\res\!S=\tau\}$.

The reader is assumed to know basics of 
nonstandard analysis and be familiar with nonstandard universes
and Loeb space construction. We suggest the reader consult
[L] and [SB] for information on those subjects. The reader
is also assumed to have basic knowledge on set theory and forcing.
The reader is recommended to consult [K1] for that.

\section{Towards Compactness}

We would like to list three facts about hyperfinite Loeb spaces 
$(\Omega,\Sigma,L_P)$, which will be used frequently throughout
this paper. 

\medskip

{\bf Fact 1}: For any $S\in\Sigma$ and any $\epsilon>0$ there
exists an $A\in\; ^*{\cal P}(\Omega)$ such that $A\subseteq S$
and $L_P(S\smallsetminus A)<\epsilon$.

\medskip

{\bf Fact 2}: For any internal sets $A_n\subseteq\Omega$ 
there exists an internal set $B\subseteq\bigcap_{n\in\omega}A_n$ 
such that $L_P(\bigcap_{n\in\omega}A_n)=L_P(B)$.

\medskip

A $\lambda$-sequence $\langle A_{\alpha}:\alpha\in\lambda\rangle$
of measurable subsets of $\Omega$ is called independent  
if for any finite $I_0\subseteq\lambda$
\[L_P(\bigcap_{\alpha\in I_0}A_{\alpha})=
\prod_{\alpha\in I_0}L_P(A_{\alpha}).\]

{\bf Fact 3}: Suppose $|\Omega|=\lambda$. Then there exists
an independent $\lambda$-sequence of internal sets of measure $\frac{1}{2}$
on $\Omega$.

\medskip

Fact 1 and Fact 2 are direct consequences of $\omega_1$-saturation
and Loeb measure construction. Fact 3 can be proved by 
finite combinatorics and the transfer property. For a Loeb space
in the general sense Fact 1 and Fact 2 are also true. But Fact 3 may
not hold. So whether or not a result about hyperfinite Loeb space
can be generalized to a general Loeb space may depend on the truth of
Fact 3.

\medskip

A set $t\subseteq \;^{<\omega}2$ is called a tree if for any
$s,s'\in\; ^{<\omega}2$, $s\subseteq s'$ and $s'\in t$ imply
$s\in t$. We use capital letter 
$T\subseteq\; ^{<\omega}2$ exclusively for
a tree with no maximal node. So every branch of $T$ is infinite.
For a tree $T$ we write $[T]$ for the set of all its branches.
In fact, every closed subset of $^{\omega}2$ could be written
as $[T]$ for some tree $T$.

\begin{definition}

A sequence of trees $\langle T_{\alpha,n}:\alpha\in\lambda\wedge
n\in\omega\rangle$ is called a $(\kappa,\lambda)$-witness if

(1) $\nu_{\omega}([T_{\alpha,n}])>\frac{n}{n+1}$,

(2) $(\forall f\in\; ^{\omega}2) (|\{\alpha\in\lambda:
\exists n (f\in [T_{\alpha,n}])\}|<\kappa)$.

\end{definition}

\begin{theorem}

Suppose there exists a $(\kappa,\lambda)$-witness
$\langle T_{\alpha,n}:\alpha\in\lambda\wedge
n\in\omega\rangle$ for some uncountable cardinals 
$\kappa$ and $\lambda$. Suppose $^*V$ is $\kappa$-saturated
and $|^*{\cal P}(\Omega)|=\lambda$. Then $\Omega$
is compact.

\end{theorem}

\noindent {\bf Proof:}\quad
Choose $\langle A_n:n\in\omega\rangle$, an independent 
$\omega$-sequence of internal subsets with measure $\frac{1}{2}$
on $\Omega$. We write $A^0_n=A_n$ and $A^1_n=\Omega\smallsetminus A_n$.
Then for any finite $s\subseteq\omega$ and any $h\in\; ^s2$
we have \[L_P(\bigcap_{n\in s}A^{h(n)}_n)=2^{-|s|}.\]
For any tree $T$ define
\[A_T=\bigcap_{n\in\omega}\bigcup_{\eta\in \;^n2\cap T}
\bigcap_{i=0}^{n-1}A^{\eta(i)}_i.\]
It is easy to see that $L_P(A_T)=\nu_{\omega}([T])$.
Note that $A_T$ is a countable intersection of internal sets.
We now want to construct an inner-regular compact family $\cal C$
of internal subsets on $\Omega$. 
Let $^*{\cal P}(\Omega)=\{a_{\alpha}:\alpha\in\lambda\}$.
For each $\alpha\in\lambda$ and $n\in\omega$ let 
$b_{\alpha,n}\subseteq a_{\alpha}\cap A_{T_{\alpha,n}}$
be internal such that \[L_P(b_{\alpha,n})=
L_P(a_{\alpha}\cap A_{T_{\alpha,n}}).\]
Then let
\[{\cal C}=\{A^l_n:n\in\omega\wedge l=0,1\}\cup\{b_{\alpha,n}:
\alpha\in\lambda\wedge n\in\omega\}.\]

\medskip

{\bf Claim 2.1}:\quad $\cal C$ is an inner-regular compact family
on $\Omega$.

Proof of Claim 2.1:\quad
The inner-regularity is clear. We need to prove the compactness.
Let ${\cal D}\subseteq {\cal C}$ be such that $\cal D$ has f.i.p.
We want to show that $\bigcap{\cal D}\not=\emptyset$.
Without loss of generality we assume that $\cal D$ is maximal.
So for each $n\in\omega$ either $A^0_n\in {\cal D}$ or
$A^1_n\in {\cal D}$ but not both. Let $h\in\; ^{\omega}2$ be such that
for each $n\in\omega$ we have $A^{h(n)}_n\in {\cal D}$.
Given any $b_{\alpha,n}\in {\cal D}$,
we want to show that $h\in [T_{\alpha,n}]$. Let $k\in\omega$.
Then \[(\bigcap_{i=0}^{k-1}A^{h(i)}_i)\cap b_{\alpha,n}\not=\emptyset.\]
So we have
\[(\bigcap_{i=0}^{k-1}A^{h(i)}_i))\cap (\bigcap_{m\in\omega}
\bigcup_{\eta\in\; ^m2\cap T_{\alpha,n}}\bigcap_{i=0}^{m-1}
A^{\eta(i)}_i)\not=\emptyset.\]
This implies that there exists an $\eta\in T_{\alpha,n}\cap\; ^k2$
such that
\[(\bigcap_{i=0}^{k-1}A^{h(i)}_i)\cap (\bigcap_{i=0}^{k-1}
A^{\eta(i)}_i)\not=\emptyset.\]
Hence we have $h\!\res\!k=\eta\!\res\!k\in T_{\alpha,n}$.
This is true for any $k\in\omega$. So $h\in [T_{\alpha,n}]$.
But we assumed that \[|\{\alpha:\exists n
(h\in [T_{\alpha,n}])\}|<\kappa.\] So $|{\cal D}|<\kappa$.
Now using $\kappa$-saturation, we get $\bigcap {\cal D}\not=\emptyset$.
\quad $\Box$

\medskip

\noindent {\bf Remark:}\quad From the definition of the compactness
we do not have to choose $\cal C$ as a family of internal sets. 
We do that because internal sets are more interesting. In this paper
if we construct a compact family we always construct a family of 
internal sets.

\begin{corollary}

Suppose CH holds and $|^*{\Bbb N}|=\omega_1$. Then
every hyperfinite Loeb space in $^*V$ is compact.

\end{corollary}

\noindent {\bf Proof:}\quad
It suffices to construct an $(\omega_1,\omega_1)$-witness.
Let $^{\omega}2=\{f_{\alpha}:\alpha\in\omega_1\}$. For each
$n\in\omega$ and $\alpha\in\omega_1$, choose $T_{\alpha,n}$
such that \[\nu_{\omega}([T_{\alpha,n}])>\frac{n}{n+1}\] and
\[[T_{\alpha,n}]\cap\{f_{\beta}:\beta\in\alpha\}=\emptyset.\]
It is clear that $\langle T_{\alpha,n}:\alpha\in\omega_1\wedge
n\in\omega\rangle$ is an $(\omega_1,\omega_1)$-witness.
\quad $\Box$

\medskip

\noindent {\bf Remarks:}\quad (1) The condition 
$|^*{\Bbb N}|=\omega_1$ implies CH by $\omega_1$-saturation
of $^*V$.

(2) If $^*V$ is an $\omega$-ultrapower of the standard universe,
then $|^*{\Bbb N}|=\omega_1$, provided CH holds.

\begin{corollary}

Suppose for any $S\subseteq\; ^{\omega}2$, $|S|<2^{\omega}$ implies
$\nu_{\omega}(S)=0$. Suppose $^*V$ is $2^{\omega}$-saturated and
$|^*{\Bbb N}|=2^{\omega}$. Then every hyperfinite 
Loeb space in $^*V$ is compact.

\end{corollary}

\noindent {\bf Proof:}\quad
By same construction as in Corollary 3 with length $2^{\omega}$
we can find a $(2^{\omega},2^{\omega})$-witness.
Now the corollary follows from $2^{\omega}$-saturation of $^*V$
and Theorem 2. \quad $\Box$

\medskip

\noindent {\bf Remark:}\quad MA implies $\nu_{\omega}(S)=0$
for any set $S\subseteq\;^{\omega}2$ with $|S|<2^{\omega}$. MA implies
also $2^{\kappa}=2^{\omega}$ for any $\kappa<2^{\omega}$, which
guarantees the existence of $2^{\omega}$-saturated nonstandard
universes.

\begin{corollary}

Suppose $\cal M$ is obtained by adding
$\lambda$ Cohen reals to a ZFC model $\cal N$ for some
$\lambda\geqslant (2^{\omega})^{\cal N}$ 
with $\lambda^{\omega}=\lambda$ in $\cal N$. 
Suppose $|^*{\Bbb N}|=\lambda$ ($2^{\omega}=\lambda$ in $\cal M$). 
Then every hyperfinite Loeb space in $^*V$ is compact.

\end{corollary}

\noindent {\bf Proof:}\quad
It suffices to construct an $(\omega_1,\lambda)$-witness.
Work in $\cal N$. For each $n\in\omega$ let 
\[{\cal T}_n=\{t\subseteq\; ^{<\omega}2:(\exists T\subseteq\; ^{<\omega}2)
(\nu_{\omega}([T])>\frac{n+1}{n+2}\wedge\exists m
(t=T\!\res\!m))\}\]
be a forcing notion ordered by the reverse of end-extension of trees
(we assume smaller conditions are stronger).
It is clear that ${\cal T}_n$ is countable and separative.
So forcing with ${\cal T}_n$ is same as adding a Cohen real.
Let ${\cal T}^{\alpha}_n={\cal T}_n$, let ${\Bbb P}_{\alpha}=
\prod_{n\in\omega}{\cal T}_n^{\alpha}$ with finite supports for each
$\alpha\in\lambda$ and let 
${\Bbb P}=\prod_{\alpha\in\lambda}{\Bbb P}_{\alpha}$ with
finite supports.
Without loss of generality we assume that ${\cal M}={\cal N}[G]$,
where $G\subseteq {\Bbb P}$ is an $\cal N$-generic filter.
For each $\alpha\in\lambda$ and $n\in\omega$ let
\[T_{\alpha,n}=\bigcup (G\cap {\cal T}^{\alpha}_n).\]
We want to show that the sequence $\langle T_{\alpha,n}
:\alpha\in\lambda\wedge n\in\omega\}$ is an
$(\omega_1,\lambda)$-witness.

Given any $f\in\; ^{\omega}2$ in $\cal M$, there exists a countable
set $S\subseteq\lambda$ in $\cal N$ such that $f\in {\cal N}[G_S]$,
where $G_S=G\cap(\prod_{\alpha\in S}{\Bbb P}_{\alpha})$.
For any $\alpha\in\lambda\smallsetminus S$ and any $n\in\omega$,
$G\cap {\cal T}^{\alpha}_n\subseteq {\cal T}^{\alpha}_n$ is a 
${\cal N}[G_S]$-generic filter. Define
\[D_f=\{t\in {\cal T}^{\alpha}_n:\exists m\;
(t \mbox{ has height }m\wedge f\!\res\!m\not\in t)\}.\]

\medskip

{\bf Claim 5.1}\quad ${\cal D}_f$ is dense in ${\cal T}^{\alpha}_n$.

Proof of Claim 5.1:\quad
Let $t\in {\cal T}^{\alpha}_n$. We want to find a 
$t'\in {\cal T}^{\alpha}_n\cap D_f$ such that $t'$ is an end-extension 
of $t$. Let $m'$ be the height of $t$. Without loss of generality
we assume that $f\!\res\!m'\in t$.
Let \[T=\{\eta\in\; ^{<\omega}2:\eta\in t\vee\eta\!\res\!m'\in t\}.\]
It is clear that $\nu_{\omega}([T])>\frac{n+1}{n+2}$.
Let \[\epsilon=\nu_{\omega}([T])-\frac{n+1}{n+2}\] and let $n'>m'$
be large enough so that $2^{-n'}<\epsilon$. Let 
\[T'=\{\eta\in T:
|\eta|\leqslant n'\vee\eta\!\res\!(n'+1)\not=f\!\res\!(n'+1)\}.\]
Now \[\nu_{\omega}([T'])>\nu_{\omega}([T])-\epsilon=\frac{n+1}{n+2}.\]
Let $m=n'+1$. Then we have 
\[t'=T'\!\res\!m\in {\cal T}^{\alpha}_n\cap D_f\]
and that $t'$ is an end-extension of $t$.
\quad $\Box$ (Claim 5.1)

\medskip

Since $D_f$ is dense in ${\cal T}^{\alpha}_n$, then
$G\cap {\cal T}^{\alpha}_n\cap D_f\not=\emptyset$. This implies
$f\not\in [T_{\alpha,n}]$. So
\[\{\alpha\in\lambda:\exists n (f\in [T_{\alpha,n}])\}
\subseteq S.\] This shows that
$\langle T_{\alpha,n}:\alpha\in\lambda\wedge n\in\omega\rangle$
is an $(\omega_1,\lambda)$-witness.\quad $\Box$

\medskip

\noindent {\bf Remarks:}\quad (1) We don't have requirements
for $^*V$ because $^*V$ is always $\omega_1$-saturated.

(2) Above three corollaries can be easily generalized to
general Loeb spaces as long as their cardinalities are $2^{\omega}$.
For example, above three corollaries are also true if we replace
a hyperfinite Loeb space by a Loeb space generated by
a nonstandard version of Lebesgue measure on unit interval. From now
on we will not make similar remarks like this. The reader should
be able to do so by himself.

\medskip

Next we will mention a property of nonstandard universes
called the $\aleph_0$-special model axiom (see [R5] or [J] for details).
The special model axiom is in fact an axiomatization of special
nonstandard universes as models (see [CK] for the definition of
a special model).
In the proof we need only some simple consequences of the property.
Let's list all the consequences we need. If $^*V$ satisfies 
the $\aleph_0$-special model axiom, then 

(1) all infinite internal sets have same cardinality, say $\lambda$,

(2) for every hyperfinite internal space 
$(\Omega,\;^*{\cal P}(\Omega),P)$ there exists a sequence
$\langle(\Omega_{\alpha},\Sigma_{\alpha}):\alpha\in cf(\lambda)\rangle$,
called a specializing sequence,
such that 

\par\par (a) $\Omega=\bigcup_{\alpha\in cf(\lambda)}\Omega_{\alpha}$,
 
\par\par (b) $^*{\cal P}(\Omega)=
\bigcup_{\alpha\in cf(\lambda)}\Sigma_{\alpha}$,

\par\par (c) if $\{A_n:n\in\omega\}\subseteq\Sigma_{\alpha}$, 
then there exists
a $B\in\Sigma_{\alpha+1}$ such that $B\subseteq\bigcap_{n\in\omega}A_n$
and $L_P(B)=L_P(\bigcap_{n\in\omega}A_n)$, 

\par\par (d) if ${\cal D}\subseteq
\Sigma_{\alpha}$ and ${\cal D}$ has f.i.p., then 
$(\bigcap {\cal D})\cap\Omega_{\alpha+1}\not=\emptyset$.

\begin{theorem}

Suppose there exists an increasing sequence 
$\langle Z_{\alpha}\subseteq\; ^{\omega}2:\alpha\in\kappa\rangle$ for
some regular cardinal $\kappa>\omega$ 
such that $\nu_{\omega}(Z_{\alpha})=0$
for every $\alpha\in\kappa$ and $^{\omega}2=
\bigcup_{\alpha\in\kappa}Z_{\alpha}$ (so $\kappa\leqslant 2^{\omega}$). 
Suppose $\lambda$ is a strong limit cardinal with $cf(\lambda)=\kappa$.
Suppose $^*V$ satisfies the $\aleph_0$-special model axiom and
$|^*{\Bbb N}|=\lambda$. Then every hyperfinite Loeb space in $^*V$ is compact.

\end{theorem}

\noindent {\bf Proof:}\quad
Given a hyperfinite Loeb space $\Omega$, let 
$\langle (\Omega_{\beta},\Sigma_{\beta}):\beta\in\kappa\rangle$ be
a specializing sequence of $\Omega$. For each $\beta\in\kappa$
let $T_{\beta,n}\subseteq\; ^{<\omega}2$ be such that 
\[\nu_{\omega}([T_{\beta,n}])>\frac{n}{n+1}\mbox{ and }
[T_{\beta,n}]\cap Z_{\beta}=\emptyset.\] 
Without loss of generality we can pick an independent
sequence $\langle A_n:n\in\omega\rangle$ of internal sets with
measure $\frac{1}{2}$ in $\Sigma_0$.
Let $^*{\cal P}(\Omega)=\{a_{\alpha}:\alpha\in\lambda\}$.
For each $\alpha\in\lambda$ let
\[g(\alpha)=\min\{\beta\in\kappa:a_{\alpha}\in\Sigma_{\beta}\}.\]
We now construct an inner-regular compact family $\cal C$ on $\Omega$.
For each $\alpha\in\lambda$ and each $n\in\omega$ let 
$b_{\alpha,n}\in\Sigma_{g(\alpha)+1}$ be such that \[b_{\alpha,n}
\subseteq a_{\alpha}\cap A_{T_{g(\alpha),n}}\] and
\[L_P(b_{\alpha,n})=L_P(a_{\alpha}\cap A_{T_{g(\alpha),n}})\]
(check Theorem 2 for the definition of $A_T$).
Define now
\[{\cal C}=\{b_{\alpha,n}:\alpha\in\lambda\wedge n\in\omega\}
\cup\{A^l_n:n\in\omega\wedge l=0,1\}\]
(recall $A^0=A$ and $A^1=\Omega\smallsetminus A$).

\medskip

{\bf Claim 6.1}\quad
$\cal C$ is an inner-regular compact family on $\Omega$.

Proof of Claim 6.1:\quad Again the inner-regularity is clear.
Let ${\cal D}\subseteq {\cal C}$ be such that $\cal D$ has
f.i.p. We want to show $\bigcap {\cal D}\not=\emptyset$.
Again we assume $\cal D$ is maximal. Let 
\[\delta=\bigcup\{g(\alpha):\exists n
(b_{\alpha,n}\in {\cal D})\}.\]

Case 1: \quad $\delta<\kappa$. \quad
Then ${\cal D}\subseteq\Sigma_{\delta+1}$. By the special model axiom
we have $\bigcap {\cal D}\not=\emptyset$.

\medskip

Case 2: \quad $\delta=\kappa$.\quad
Let $h\in\; ^{\omega}2$ be such that $A^{h(n)}_n\in {\cal D}$.
Same as the proof of Claim 2.1 we have that $h\in [T_{g(\alpha),n}]$
if $b_{\alpha,n}\in {\cal D}$. But there is a $\beta_0\in\kappa$
such that $h\in Z_{\beta_0}$. So we have $h\not\in
[T_{g(\alpha),n}]$ for any $g(\alpha)>\beta_0$. This contradicts
$\delta=\kappa$.
\quad $\Box$ 

\medskip

\noindent {\bf Remarks:}\quad 
(1) The nonstandard universes satisfying the $\aleph_0$-special
model axiom and having cardinality $\lambda$ 
exist. In fact, those universes are 
frequently used by Ross (see [R4] and [R5]).

(2) This theorem guarantees the existence of arbitrarily large
compact hyperfinite Loeb spaces.

(3) The set theoretical assumption besides ZFC for $\cal M$ 
in this theorem is rather weak. The model $\cal M$ satisfies 
the assumption if $\cal M$ is a model of {\em e.g.} CH or MA,
or is obtained by adding enough Cohen reals.

\begin{corollary}

Suppose $\cal M$ is same as in Corollary 5.
In $\cal M$ suppose $\lambda$ is a strong limit cardinal 
with $cf(\lambda)\leqslant\kappa$.
Suppose $|^*{\Bbb N}|=\lambda$ and $^*V$ satisfies the 
$\aleph_0$-special model axiom. 
Then every hyperfinite Loeb space in $^*V$ is compact.

\end{corollary}

\noindent {\bf Proof:}\quad
First we arrange the Cohen forcing such that $\cal M$ is a forcing
extension of some model of ZFC by adding $cf(\lambda)$ Cohen reals.
Then by [K2, Theorem 3.20] we know that $^{\omega}2$ is a union of
an increasing $cf(\lambda)$-sequence of measure zero sets.
\quad $\Box$

\medskip

\noindent {\bf Remark:}\quad There is another proof by a method
similar to the proof of Corollary 5.

\begin{theorem}

Suppose $\kappa$ is a strong limit cardinal with $cf(\kappa)=\omega$
and suppose $\lambda=\kappa^+=2^{\kappa}$. Suppose $^*V$
is $\lambda$-saturated and $|^*{\Bbb N}|=\lambda$. Then every
hyperfinite Loeb space in $^*V$ is compact.

\end{theorem}

\noindent {\bf Proof:}\quad Given a hyperfinite Loeb space $\Omega$
in $^*V$, we want to show that $\Omega$ is compact.
Let $\kappa=\bigcup_{n\in\omega}\kappa_n$ be such that
$\kappa_0>\omega$ and $\kappa_{n+1}>2^{\kappa_n}$ for each $n\in\omega$.
Choose an independent $\kappa$-sequence 
\[\langle A_{\alpha}\in\;^*{\cal P}(\Omega):\alpha\in\kappa\rangle\]
of internal sets with measure $\frac{1}{2}$. For any $n,m$ let
\[{\cal B}_n={\cal B}(\{A_{\alpha}:\alpha\in\kappa_{n+1}\})\] be
the Boolean algebra generated by $A_{\alpha}$'s for all 
$\alpha\in\kappa_{n+1}$ and let
\[Pos_{n,m}=\{X\in {\cal B}_n: L_P(X)>\frac{m}{m+1}\}.\]
Note that every $X\in Pos_{n,m}$ is internal with measure $>\frac{m}{m+1}$
because it is a finite Boolean combination of internal sets.
For each $n,m\in\omega$ let
\[I_{n,m}=\{E\subseteq Pos_{n,m}:E\mbox{ has f.i.p.}\}.\]
For each $m\in\omega$ let
\[Pos_{<\omega,m}=\{\bar{X}:\exists n (\bar{X}=\langle
X_0,X_1,\ldots,X_{n-1}\rangle\wedge 
(\forall i<n)(X_i\in Pos_{n,m}))\}.\]
and let \[{\cal F}_m=\{F:F\mbox{ is a function from }Pos_{<\omega,m}
\mbox{ to }\bigcup_{n\in\omega}I_{n,m}\mbox{ such that }\]
\[(\forall\bar{X}=\langle X_0,\ldots,X_{n-1}\rangle\in Pos_{<\omega,m})
(F(\bar{X})\in I_{n,m})\}.\]
It is clear that $|{\cal F}_m|\leqslant\kappa^\kappa=\lambda$.
Let ${\cal F}_m=\{F_{\alpha,m}:\alpha\in\lambda\}$ be a fixed
enumeration. For each $\alpha\in\lambda$ let's fix an increasing
sequence $\langle B_{\alpha,n}\subseteq\alpha:n\in\omega\rangle$
such that $|B_{\alpha,n}|\leqslant\kappa_n$ for each $n\in\omega$
and $\alpha=\bigcup_{n\in\omega}B_{\alpha,n}$.
We define a function $f_{\alpha,m}$ from $\omega$ to 
$\bigcup_{n\in\omega}Pos_{n,m}$ for each $\alpha\in\lambda$
by induction on $n$ such that for each $n\in\omega$

(1) $f_{\alpha,m}(n)\in Pos_{n,m}$,

(2) $f_{\alpha,m}(n+1)\subseteq f_{\alpha,m}(n)$,

(3) $f_{\alpha,m}(n)\not\in\bigcup\{F_{\beta,m}
(f_{\alpha,m}\!\res\!n):\beta\in B_{\alpha,n}\}$.

Suppose we have defined $f_{\alpha,m}\!\res\!n$.

\medskip

{\bf Claim 8.1}\quad There is an $X=f_{\alpha,m}(n)$ such that
(1), (2) and (3) hold.

Proof of Claim 8.1:\quad For each $\beta\in B_{\alpha,n}$ let
${\cal D}_{\beta}$ be an ultrafilter on ${\cal B}_n$ such that 
\[{\cal D}_{\beta}\supseteq F_{\beta,m}(f_{\alpha,m}\!\res\!n)\]
(note that $F_{\beta,m}(f_{\alpha,m}\!\res\!n)$ has f.i.p.).
Let $f_{\alpha,m}(n-1)=C$ at stage $n>0$ (replace $C$ by $\Omega$
at stage $n=0$).
For each $\gamma\in [\kappa_n,\kappa_{n+1})$ let
\[J_{\gamma}=\{\beta\in B_{\alpha,n}:
C\cap A_{\gamma}\in {\cal D}_{\beta}\}.\]
Since $|B_{\alpha,n}|\leqslant\kappa_n$ and 
$2^{\kappa_n}<\kappa_{n+1}=|[\kappa_n,\kappa_{n+1})|$, then there
exists an $E\subseteq [\kappa_n,\kappa_{n+1})$ with $|E|=\kappa_{n+1}$
such that for any two different $\gamma,\gamma'\in E$ we have
$J_{\gamma}=J_{\gamma'}$.
Let $\gamma_0<\gamma'_0<\gamma_1<\gamma'_1<\ldots$ be in $E$ and let
\[C_n=(A_{\gamma_n}\cap C)\Delta (A_{\gamma'_n}\cap C)=
(A_{\gamma_n}\Delta A_{\gamma'_n})\cap C,\]
where $\Delta$ means symmetric difference.
It is easy to see that for any $n\in\omega$ and any 
$\beta\in B_{\alpha,n}$ we have $C_n\not\in {\cal D}_{\beta}$.
It is also clear that 
\[L_P(C_n)=L_P(C)L_P(A_{\gamma_n}\Delta A_{\gamma'_n})
=\frac{1}{2}L_P(C).\]
Since $L_P(C)>\frac{m}{m+1}$, there exists a big enough $N\in\omega$
such that
\[(1-(\frac{1}{2})^N)L_P(C)>\frac{m}{m+1}.\]
Now let $f_{\alpha,m}(n)=\bigcup_{i=0}^{N-1} C_i$.
It is easy to see that (2) and (3) hold.
For (1) we have
\[L_P(f_{\alpha,m}(n))=L_P(C\cap (\bigcup_{i=0}^{N-1}
(A_{\gamma_i}\Delta A_{\gamma'_i})))\]
\[=L_P(C)L_P(\bigcup_{i=0}^{N-1}(A_{\gamma_i}\Delta A_{\gamma'_i}))\]
\[=(1-(\frac{1}{2})^N)L_P(C)>\frac{m}{m+1}.\]
\quad $\Box$ (Claim 8.1)

\medskip

We now define an inner-regular compact family $\cal C$ on $\Omega$.
Let $^*{\cal P}(\Omega)=\{a_{\alpha}:\alpha\in\lambda\}$ be
an enumeration. For each $\alpha\in\lambda$ and $m\in\omega$, let
\[b_{\alpha,m}\subseteq a_{\alpha}\cap 
(\bigcap_{n\in\omega}f_{\alpha,m}(n))\] be internal such that
\[L_P(b_{\alpha,m})=L_P(a_{\alpha}\cap 
(\bigcap_{n\in\omega}f_{\alpha,m}(n))).\] Now let
\[{\cal C}=\{b_{\alpha,m}:\alpha\in\lambda\wedge m\in\omega\}.\]

{\bf Claim 8.2}\quad
$\cal C$ is an inner-regular compact family on $\Omega$.

Proof of Claim 8.2:\quad
Let ${\cal D}\subseteq {\cal C}$ be such that $\cal D$ has f.i.p.
If $|{\cal D}|<\lambda$, then $\bigcap {\cal D}\not=\emptyset$
by $\lambda$-saturation. So let's assume $|{\cal D}|=\lambda$.
Hence there exists an $m_0\in\omega$ such that
\[Z=\{\alpha:b_{\alpha,m_0}\in {\cal D}\}\] has cardinality $\lambda$.
\quad $\Box$ (Claim 8.2)
\medskip

{\bf Claim 8.3}\quad
There exists an $\bar{X}=\langle X_0,X_1,\ldots,X_{n-1}\rangle$
for some $n\in\omega$ such that
\[\{f_{\alpha,m_0}(n): \alpha\in Z\wedge 
f_{\alpha,m_0}\!\res\!n=\bar{X}\}\not\in I_{n,m_0}.\]

Proof of Claim 8.3:\quad
Suppose not. Then we can define a function
\[F:Pos_{<\omega,m_0}\mapsto\bigcup_{n\in\omega}I_{n,m_0}\] such that
for each $\bar{X}=\langle X_0,\ldots,X_{n-1}\rangle$
\[F(\bar{X})=\{f_{\alpha,m_0}(n):\alpha\in Z\wedge
f_{\alpha,m_0}\!\res\!n=\bar{X}\}.\]
It is clear that $F\in {\cal F}_{m_0}$. So there is an 
$\beta\in\lambda$ such that $F=F_{\beta,m_0}$.
Since $|Z|=\lambda$ there is an $\alpha\in Z$ such that
$\alpha>\beta$. Now choose large enough $n\in\omega$ such that
$\beta\in B_{\alpha,n}$. Then \[f_{\alpha,m_0}(n)\not\in
F_{\beta,m_0}(f_{\alpha,m_0}\!\res\!n)\] by the construction of
$f_{\alpha,m_0}$. But this contradicts the definition
of $F=F_{\beta,m_0}$.\quad $\Box$ (Claim 8.3)

\medskip

We continue the proof of Claim 8.2.
By Claim 8.3 there exists an \hfill\break 
$\bar{X}=\langle X_0,\ldots,X_{n-1}\rangle$
such that \[\{f_{\alpha,m_0}(n):\alpha\in Z\wedge 
f_{\alpha,m_0}\!\res\!n=\bar{X}\}\] does not have f.i.p.
Hence $\{b_{\alpha,m_0}:\alpha\in Z\}$ 
does not have f.i.p. because \hfill\break
$b_{\alpha,m_0}\subseteq f_{\alpha,m_0}(n)$.
This contradicts that $\cal D$ has f.i.p. \quad $\Box$

\medskip

\noindent {\bf Remarks:}\quad
(1) $\lambda$-saturated nonstandard universes of cardinality $\lambda$
exist because $\lambda^{<\lambda}=\lambda$. 

(2) Under certain assumptions for $\cal M$, {\em e.g.} Singular
Cardinal Hypothesis, 
this theorem guarantees the existence of arbitrarily large
compact hyperfinite Loeb spaces with regular cardinality.

(3) The proof of this theorem is implicitly included in [Sh575].

\section{Towards non-compactness}

For a measure space $(X,\Sigma,P)$ we write $\bar{\cal B}(X)$
for the measure algebra of $X$, 
{\em i.e.} the Boolean algebra of measurable
sets modulo the ideal of measure zero sets. 
If ${\cal D}\subseteq\Sigma$, we write
$\bar{\cal B}(\cal D)$ for the complete subalgebra of 
$\bar{\cal B}(X)$ generated by $\cal D$.

\begin{theorem}

Suppose $\lambda$ is a regular cardinal such that $\kappa^{\omega}<
\lambda$ for every $\kappa<\lambda$. Suppose $|\Omega|=\lambda$.
Then $\Omega$ is not compact.

\end{theorem}

\noindent {\bf Proof:}\quad
Suppose not. Let $\cal C$ be the inner-regular compact family
on $\Omega$. Let \hfill\break
$\langle A_{\alpha}:\alpha\in\lambda\rangle$
be an independent $\lambda$-sequence of internal subsets of $\Omega$
of Loeb measure $\frac{1}{2}$. 
Pick $X^l_{\alpha,n}\in {\cal C}$ for every $\alpha\in\lambda$,
$n\in\omega$ and $l=0,1$ such that 
$X^l_{\alpha,n}\subseteq A^l_{\alpha}$ and \[L_P(\bigcup_{n\in\omega}
X^l_{\alpha,n})=\frac{1}{2}.\] Note that $A^0_{\alpha}=A_{\alpha}$ and
$A^1_{\alpha}=\Omega\smallsetminus A_{\alpha}$, but 
generally $X^1_{\alpha,n}$
will not be $\Omega\smallsetminus X^0_{\alpha,n}$. 

\medskip

{\bf Claim 9.1}\quad There exists an $E\in [\lambda]^{\lambda}$ and
there exists an $n(\alpha,l)$ for each $\alpha\in E$ and $l=0,1$
such that for any $m\in\omega$, any distinct
$\{\alpha_0,\ldots,\alpha_{m-1}\}\subseteq E$ and any $h\in\; ^m2$ we have
\[L_P(\bigcap_{i=0}^{m-1}X^{h(i)}_{\alpha_i,n(\alpha_i,h(i))})>0.\]

\medskip

The theorem follows from the claim. Since $|\Omega|=\lambda<2^{\lambda}$,
we can find  $f\in\; ^E2$ such that \[\bigcap_{\alpha\in E}
X^{f(\alpha)}_{\alpha,n(\alpha,f(\alpha))}=\emptyset.\]
But $\{X^{f(\alpha)}_{\alpha,n(\alpha,f(\alpha))}:\alpha\in E\}\subseteq
{\cal C}$ has f.i.p.

\medskip

Proof of Claim 9.1:\quad For each measurable set $A\subseteq\Omega$
let $\bar{A}$ denote the element in the measure algebra, represented
by $A$. For each $\alpha\in\lambda$ let
\[\bar{\cal B}_{\alpha}=\bar{\cal B}(\{A_{\beta}:
\beta<\alpha\}\cup\{X^l_{\beta,n}:
\beta<\alpha\wedge n\in\omega\wedge l=0,1\}).\]
Recall that $\bar{\cal B}(X)$ for some family $X$ of measurable sets
is a complete subalgebra of measure algebra on $\Omega$ generated
by $X$. By c.c.c. of measure algebra it is easy to see that
$|\bar{\cal B}_{\alpha}|\leqslant |\alpha|^{\omega}$.
Notice that the sequence $\langle\bar{\cal B}_{\alpha}:
\alpha\in\lambda\rangle$ is increasing and 
$\bar{\cal B}_{\alpha}=\bigcup_{\beta<\alpha}\bar{\cal B}_{\beta}$
when $cf(\alpha)>\omega$.
Let $\cal B$ be a complete Boolean algebra. Given $a\in {\cal B}$
and a subalgebra ${\cal B}'\subseteq {\cal B}$, define a function
\[g(a,{\cal B}')=\inf\{b\in {\cal B}':b\geqslant a\}.\]
$g(a,{\cal B}')$ exists since $\cal B$ is complete.
Let \[D=\{\alpha\in\lambda:cf(\alpha)>\omega\}.\] Then $D$ is stationary
in $\lambda$. For each $\alpha\in D$ let
\[d(\alpha)=\min\{\beta:g(\bar{A}^0_{\alpha},\bar{\cal B}_{\alpha})
\in\bar{\cal B}_{\beta}\wedge 
g(\bar{A}^1_{\alpha},\bar{\cal B}_{\alpha})\in\bar{\cal B}_{\beta}\}.\] 
Then $d(\alpha)<\alpha$ for every $\alpha\in D$. 
By Pressing-Down Lemma we can find a stationary subset $E\subseteq D$
and an $\alpha_0\in\lambda$ such that $d(\alpha)=\alpha_0$ for every
$\alpha\in E$. Since \[|\bar{\cal B}_{\alpha_0}|\leqslant 
|\alpha_0|^{\omega}<\lambda,\] we can assume that there are 
$b_0,b_1\in\bar{\cal B}_{\alpha_0}$ such that for all $\alpha\in E$
we have \[g(\bar{A}^0_{\alpha},\bar{\cal B}_{\alpha})=b_0\mbox{ and }
g(\bar{A}^1_{\alpha},\bar{\cal B}_{\alpha})=b_1.\]
By thinning $E$ further we can assume that $\bar{A}_{\alpha}\not\in
\bar{\cal B}_{\alpha_0}$ for each $\alpha\in E$. Hence
$\bar{A}_{\alpha}\not\in\bar{\cal B}_{\alpha}$ for each $\alpha\in E$. 
It is easy to see that $b_0\wedge b_1\not=0$ because otherwise
we have, for any $\alpha\in E$,
\[\bar{A}^0_{\alpha}\leqslant b_0\leqslant -b_1\leqslant 
-\bar{A}^1_{\alpha}=\bar{A}^0_{\alpha}\] and this implies
$\bar{A}^0_{\alpha}=b_0\in\bar{\cal B}_{\alpha_0}$.
 
\medskip

{\bf Claim 9.2}\quad
For any $\alpha\in E$ there exist $n(\alpha,0)$ and $n(\alpha,1)$
such that 
\[g(\bar{X}^0_{\alpha,n(\alpha,0)},\bar{\cal B}_{\alpha})\wedge
g(\bar{X}^1_{\alpha,n(\alpha,1)},\bar{\cal B}_{\alpha})\not=0.\]

Proof of Claim 9.2:\quad Suppose not. Then for any $n,m\in\omega$
we have \[\bar{X}^0_{\alpha,n}\leqslant
g(\bar{X}^0_{\alpha,n},\bar{\cal B}_{\alpha})\leqslant
-g(\bar{X}^1_{\alpha,m},\bar{\cal B}_{\alpha})\leqslant
-\bar{X}^1_{\alpha,m}.\]
So then
\[\bar{A}^0_{\alpha}=\bigvee_{n\in\omega}\bar{X}^0_{\alpha,n}\leqslant
\bigvee_{n\in\omega}g(\bar{X}^0_{\alpha,n},\bar{\cal B}_{\alpha})\leqslant
\bigwedge_{m\in\omega}(-g(\bar{X}^1_{\alpha,m},\bar{\cal B}_{\alpha}))\]
\[\leqslant\bigwedge_{m\in\omega}(-\bar{X}^1_{\alpha,m})=
-\bigvee_{m\in\omega}\bar{X}^1_{\alpha,m}=-\bar{A}^1_{\alpha}=
\bar{A}^0_{\alpha}.\]
This implies \[\bar{A}_{\alpha}=\bar{A}^0_{\alpha}=
\bigvee_{n\in\omega}g(\bar{X}^0_{\alpha,n},\bar{\cal B}_{\alpha})\in
\bar{\cal B}_{\alpha},\]
a contradiction. \quad $\Box$ (Claim 9.2)

\medskip

We continue to prove Claim 9.1. By thinning $E$ even further we can assume
that there exist $d_0,d_1\in\bar{\cal B}_{\alpha_0}$ such that
for every $\alpha\in E$,
\[g(\bar{X}^0_{\alpha,n(\alpha,0)},\bar{\cal B}_{\alpha})=d_0\mbox{ and }
g(\bar{X}^1_{\alpha,n(\alpha,1)},\bar{\cal B}_{\alpha})=d_1.\]
Clearly $d=d_0\wedge d_1\not=0$.

We now prove, by induction on $m$, that for any $m\in\omega$,
for any distinct \hfill\break 
$\{\alpha_0,\ldots,\alpha_{m-1}\}\subseteq E$
and for any $h\in\; ^m2$, 
\[d\wedge(\bigwedge_{i=0}^{m-1}\bar{X}^{h(i)}_{\alpha_i,
n(\alpha_i,h(i))})\not=0.\]

Suppose the above is true for $m$, but not true for $m+1$.
Pick some distinct $\{\alpha_0,\ldots,\alpha_m\}\subseteq E$
and $h\in\; ^{m+1}2$ such that $\alpha_i<\alpha_m$ for each $i<m$ and
\[d\wedge(\bigwedge_{i=0}^m
\bar{X}^{h(i)}_{\alpha_i,n(\alpha_i,h(i))})=0.\]
Without loss of generality let $h(m)=0$. Then
\[(d\wedge(\bigwedge_{i=0}^{m-1}
\bar{X}^{h(i)}_{\alpha_i,n(\alpha_i,h(i))}))\wedge
\bar{X}^0_{\alpha_m,n(\alpha_m,0)}=0.\]
Then \[\bar{X}^0_{\alpha_m,n(\alpha_m,0)}\leqslant
-(d\wedge(\bigwedge_{i=0}^{m-1}
\bar{X}^{h(i)}_{\alpha_i,n(\alpha_i,h(i))})).\]
So \[g(\bar{X}^0_{\alpha_m,n(\alpha_m,0)},\bar{\cal B}_{\alpha_m})=d_0
\leqslant -(d\wedge(\bigwedge_{i=0}^{m-1}
\bar{X}^{h(i)}_{\alpha_i,n(\alpha_i,h(i))})).\]
This implies 
\[d_0\wedge (d\wedge(\bigwedge_{i=0}^{m-1}
\bar{X}^{h(i)}_{\alpha_i,n(\alpha_i,h(i))}))=
d\wedge(\bigwedge_{i=0}^{m-1}
\bar{X}^{h(i)}_{\alpha_i,n(\alpha_i,h(i))})=0.\]
This contradicts the inductive hypothesis.\quad $\Box$

\medskip

\noindent {\bf Remarks:}\quad
(1) Not like other results so far, Theorem 9 is a consequence of
ZFC.

(2) When $\lambda=(\eta^{\theta})^+$ for some infinite cardinals
$\eta$ and $\theta$, we have $\kappa^{\omega}<\lambda$ for any
$\kappa<\lambda$. So ZFC implies the existence of arbitrarily
large non-compact hyperfinite Loeb spaces in some nonstandard
universes.

(3) The proof of this theorem is implicitly included in [Sh92].

(4) The proof works also for general Loeb spaces if they have an
independent $\lambda$-sequence of measure $\frac{1}{2}$.

\begin{theorem}

Suppose $\lambda$ is a strong limit cardinal, $\kappa=cf(\lambda)$
and $\mu^{\omega}<\kappa$ for any $\mu<\kappa$. Suppose 
$|\Omega|=\lambda$. Then $\Omega$ is not compact.

\end{theorem}

\noindent {\bf Proof:}\quad
Suppose $\Omega$ is compact and let $\cal C$ be an inner-regular
compact family on $\Omega$. Let $\langle\lambda_{\alpha}:
\alpha\in\kappa\rangle$ be an increasing sequence such that
$2^{\lambda_{\alpha}}<\lambda_{\alpha+1}$ for each $\alpha\in\kappa$
and $\lambda=\bigcup_{\alpha\in\kappa}\lambda_{\alpha}$.
Let $\Omega=\{a_{\beta}:\beta\in\lambda\}$ be an enumeration
and let $\Omega_{\alpha}=\{a_{\beta}:\beta<\lambda_{\alpha}\}$
for each $\alpha\in\kappa$.
Choose an independent $\lambda$-sequence of internal sets of measure
$\frac{1}{2}$, say $\langle A_{\beta}:\beta\in\lambda\rangle$, on $\Omega$.

\medskip

{\bf Claim 10.1}\quad There exist two different ordinals
$\gamma_{\alpha},\gamma'_{\alpha}
\in [\lambda_{\alpha},\lambda_{\alpha+1})$ such that
$A_{\gamma_{\alpha}}\Delta A_{\gamma'_{\alpha}}\cap\Omega_{\alpha}
=\emptyset$ for each $\alpha\in\kappa$.

Proof of Claim 10.1:\quad Since $2^{\lambda_{\alpha}}
<\lambda_{\alpha+1}$, then there exist two different $\gamma_{\alpha}$
and $\gamma'_{\alpha}$ in $[\lambda_{\alpha},\lambda_{\alpha+1})$ 
such that \[A_{\gamma_{\alpha}}\cap\Omega_{\alpha}=
A_{\gamma'_{\alpha}}\cap\Omega_{\alpha}.\] 
Hence we have 
$A_{\gamma_{\alpha}}\Delta A_{\gamma'_{\alpha}}\cap\Omega_{\alpha}
=\emptyset$.\quad $\Box$ (Claim 10.1)

\medskip

For any $\alpha\in\kappa$ let $B_{\alpha}=
A_{\gamma_{\alpha}}\Delta A_{\gamma'_{\alpha}}$.
It is easy to see that $\langle B_{\alpha}:\alpha\in\kappa\rangle$
is an independent sequence of measure $\frac{1}{2}$.
By inner-regularity of $\cal C$ we can find $X_{\alpha,n}\in {\cal C}$
for each $\alpha$ and each $n\in\omega$ such that
$X_{\alpha,n}\subseteq B_{\alpha}$ and \[L_P(\bigcup_{n\in\omega}
X_{\alpha,n})=\frac{1}{2}.\]
By a similar method as in the proof of Claim 9.1 we can find
an $E\in [\kappa]^{\kappa}$ and an
$n(\alpha)\in\omega$ for each $\alpha\in E$ such that
for any $m\in\omega$ and any distinct 
$\{\alpha_0,\ldots,\alpha_{m_1}\}\subseteq E$ 
we have \[\bigcap_{i=0}^{m-1}X_{\alpha_i,n(\alpha_i)}\not=\emptyset.\]
So $\{X_{\alpha,n(\alpha)}:\alpha\in E\}\subseteq {\cal C}$ has
f.i.p., but \[\bigcap_{\alpha\in E}X_{\alpha,n(\alpha)}\subseteq
\bigcap_{\alpha\in E}B_{\alpha}=\emptyset.\]
\quad $\Box$

\medskip

\noindent {\bf Remark:} \quad Theorem 10 is also a consequence
of ZFC. So ZFC implies the existence of non-compact hyperfinite Loeb spaces
of arbitrarily large singular cardinalities.

\medskip

We need Maharam Theorem for next two theorems. 
Given a complete Boolean algebra $\cal B$. For any $X\subseteq {\cal B}$
recall that $\bar{\cal B}(X)$ is the complete subalgebra generated by $X$. Let
\[\tau({\cal B})=\min\{|X|:X\subseteq {\cal B}\wedge 
{\cal B}=\bar{\cal B}(X)\}.\] For any $a\in {\cal B}\smallsetminus\{0\}$
let ${\cal B}\!\res\!a$ be the Boolean algebra $\{b\wedge a:b\in
{\cal B}\}$ with $a$ being the largest element $1$ in ${\cal B}\!\res\! a$. 
A complete Boolean algebra $\cal B$
is called homogeneous if $\tau ({\cal B})=\tau ({\cal B}\!\res\! a)$
for every $a\in {\cal B}\smallsetminus\{0\}$. The following is a version
of Maharam Theorem (see [F, pp.911 Theorem 3.5]).

\begin{quote} 

{\bf Maharam Theorem}\quad Let $\cal B$ be a homogeneous
measure algebra of a probability space 
with $\tau(\cal B)=\lambda$. Then there is a measure
preserving isomorphism $\Phi$ from $\cal B$ to ${\cal B}(^{\lambda}2)$.

\end{quote}

Let $\mu$ be a cardinal. For next two theorems we always denote, for 
each $\alpha\in\mu$, 
\[B_{\alpha}=\{f\in\; ^{\mu}2:f(\alpha)=0\}.\]
For any set $X\subseteq\; ^{\mu}2$ let $supt(X)$ denote the support of
$X$, {\em i.e.} the smallest $w\subseteq\mu$ such that for any
$f,f'\in\; ^{\mu}2$ we have $f\!\res\!w=f'\!\res\!w$ implies
$f\in X$ iff $f'\in X$. Clearly $supt(X)$ is at most countable
if $X$ is a Baire set.
For any measurable set $X$ in a measure space
we denote again $\bar{X}$ for the element in the measure 
algebra, represented by $X$.

\begin{theorem}

Suppose $\cal M$ is obtained by adding 
$\lambda$ random reals to a ZFC model $\cal N$ for some
regular $\lambda>(2^{\omega})^{\cal N}$ with $\lambda^{\omega}=\lambda$. 
Suppose $|\Omega|=\lambda$ ($\lambda=2^{\omega}$ in $\cal M$). 
Then $\Omega$ is not compact.

\end{theorem}

\noindent {\bf Proof:}\quad Let $\Omega$ be a hyperfinite Loeb space
in $^*V$ with $|\Omega|=\lambda$. 
Without loss of generality we assume that ${\cal B}(\Omega)$
is homogeneous and $\tau({\cal B}(\Omega))=\mu$ for some $\mu\geqslant
\lambda$. The reason for that is the following. It is easy to see that 
for any internal subset $a$ of $\Omega$ with positive Loeb measure
there exists an $|a|$-independent sequence of measure 
$\frac{1}{2}L_P(a)$ on ${\cal B}(\Omega)\!\res\!\bar{a}$. 
Since $|a|=\lambda$, then
$\tau({\cal B}(\Omega)\!\res\!\bar{a})\geqslant\lambda$. 
Suppose ${\cal B}(\Omega)$ is not homogeneous. Then we can choose
an internal subset $a\subseteq\Omega$ such that $L_P(a)>0$ and
$\mu=\tau({\cal B}(\Omega)\!\res\!\bar{a})$ is the smallest.
Hence $\mu\geqslant\lambda$ and ${\cal B}(\Omega)\!\res\!\bar{a}$ 
is homogeneous. Then we could replace $\Omega$ by $a$.

By Maharam Theorem let $\Phi:{\cal B}(\Omega)\cong {\cal B}(^{\mu}2)$
be the measure preserving isomorphism. For each $\alpha\in\lambda\subseteq\mu$
let $A_{\alpha}\subseteq\Omega$ be measurable such that 
\[\Phi(\bar{A}_{\alpha})=\bar{B}_{\alpha}.\]
Suppose $\Omega$ is compact and let $\cal C$ be the inner-regular
compact family on $\Omega$. Again let $A^0_{\alpha}=A_{\alpha}$, 
$A^1_{\alpha}=\Omega\smallsetminus A_{\alpha}$, $B^0_{\alpha}=B_{\alpha}$ 
and $B^1_{\alpha}=\;^{\mu}2\smallsetminus B_{\alpha}$. Then there exist
$X^l_{\alpha,n}\in {\cal C}$ such that $X^l_{\alpha,n}\subseteq
A^l_{\alpha}$ and \[L_P(\bigcup_{n\in\omega}X^l_{\alpha,n})=\frac{1}{2},\]
where $X^1_{\alpha,n}$ may not be $\Omega\smallsetminus X^0_{\alpha,n}$.
For each $\alpha\in\lambda$, $n\in\omega$ and $l=0,1$ let
$Y^l_{\alpha,n}\subseteq B^l_{\alpha}$ be a Baire set such that
\[\Phi(\bar{X}^l_{\alpha,n})=\bar{Y}^l_{\alpha,n}.\]
We want to find an $E\in [\lambda]^{\lambda}$ and an $n(\alpha,l)$
for each $\alpha\in E$ and $l=0,1$ such that
for any $m\in\omega$, any distinct $\{\alpha_0,\ldots,\alpha_{m-1}\}
\subseteq E$ and any $h\in\; ^m2$
\[\bigwedge_{i=0}^{m-1}\bar{Y}^{h(i)}_{\alpha_i,n(\alpha_i,h(i))}
\not=0.\] This is enough to prove the theorem because by Maharam's
isomorphism we have a family 
\[{\cal F}_f=\{X^{f(\alpha)}_{\alpha,n(\alpha,f(\alpha))}:\alpha\in E\}\]
with f.i.p. for every $f\in\; ^E2$. But $|\Omega|=\lambda$. So
there must be a family ${\cal F}_f$ for some $f\in\; ^E2$ with empty
intersection, which contradicts that $\cal C$ is a compact family.

Let ${\Bbb P}={\cal B}(^{\lambda}2)$ be the forcing in $\cal N$ 
for adding $\lambda$ random reals. Since ${\Bbb P}$ has c.c.c., then
for each $\alpha\in\lambda$ there exists a countable set 
$v_{\alpha}\subseteq\lambda$ in $\cal N$ such that
\[\{Y^l_{\alpha,n}:n\in\omega\wedge l=0,1\}\subseteq 
{\cal N}^{{\cal B}(^{v_{\alpha}}2)}.\] 
Work in $\cal N$. Let $\dot{u}_{\alpha}$ be a 
$\Bbb P$-name for \[\bigcup_{n\in\omega,l=0,1}
supt(Y^l_{\alpha,n}).\] 
Again since $\Bbb P$ has c.c.c. there exists a countable
$w_{\alpha}\subseteq\mu$ such that
\[\Vdash_{\Bbb P}\dot{u}_{\alpha}\subseteq w_{\alpha}.\]
Note that $\alpha\in w_{\alpha}$ for every $\alpha\in\lambda\subseteq\mu$.
Since $\lambda>2^{\omega}\geqslant\omega_1$, 
we can find a $\bar{v}\subseteq\mu$ with $|\bar{v}|<\lambda$
and an $E\in [\lambda]^{\lambda}$ such that for any two different
$\alpha,\beta\in E$
\[w_{\alpha}\cap w_{\beta}\subseteq\bar{v}.\]
So in $\cal M$ for any two different $\alpha,\beta\in E$ we have 
\[supt(Y^l_{\alpha,n})\cap supt(Y^{l'}_{\beta,m})\subseteq\bar{v}.\]
Without loss of generality we assume that 
$\bar{v}<\lambda$ is an ordinal, $\bar{v}\cap E=\emptyset$
and \[|\lambda\smallsetminus\bigcup_{\alpha\in E}v_{\alpha}|
=\lambda.\]
For any $X\subseteq\; ^{\mu}2$ and $\eta\in\; ^{\bar{v}}2$ let
\[X(\eta)=\{\xi\in\;^{\mu\smallsetminus\bar{v}}2:\eta\hat{\;}\xi\in X\}.\]
Now we work in $\cal M$. For each $\alpha\in E$, $n\in\omega$ and $l=0,1$
let \[C^l_{\alpha,n}=\{\eta\in\; ^{\bar{v}}2:\nu_{\mu\smallsetminus
\bar{v}}(Y^l_{\alpha,n}(\eta))>0\}.\]

{\bf Claim 11.1}\quad
$\nu_{\bar{v}}(\bigcup_{n\in\omega}C^l_{\alpha,n})=1$ for each $\alpha
\in E$ and $l=0,1$.

Proof of Claim 11.1:\quad For any $\eta\in\; ^{\bar{v}}2$ we have
\[\nu_{\mu\smallsetminus\bar{v}}(\bigcup_{n\in\omega}Y^l_{\alpha,n}
(\eta))=\nu_{\mu\smallsetminus\bar{v}}
((\bigcup_{n\in\omega}Y^l_{\alpha,n})(\eta))\leqslant\frac{1}{2}\] because
$(\bigcup_{n\in\omega}Y^l_{\alpha,n})(\eta)\subseteq B^l_{\alpha}(\eta)$
and $\nu_{\mu\smallsetminus\bar{v}}(B^l_{\alpha}(\eta))=\frac{1}{2}$.
So by Fubini Theorem we have
\[\frac{1}{2}=\nu_{\mu}(\bigcup_{n\in\omega}Y^l_{\alpha,n})\]
\[=\int_{^{\bar{v}}2}\nu_{\mu\smallsetminus\bar{v}}(\bigcup_{n\in\omega}
Y^l_{\alpha,n}(\eta))d\nu_{\bar{v}}(\eta)\]
\[\leqslant\frac{1}{2}\nu_{\bar{v}}(\{\eta:\nu_{\mu\smallsetminus\bar{v}}
(\bigcup_{n\in\omega}Y^l_{\alpha,n}(\eta))>0\}).\]
This implies \[\nu_{\bar{v}}(\{\eta:\nu_{\mu\smallsetminus\bar{v}}
(\bigcup_{n\in\omega}Y^l_{\alpha,n}(\eta))>0\})=1.\] But
\[\bigcup_{n\in\omega}C^l_{\alpha,n}=\{\eta:\nu_{\mu\smallsetminus\bar{v}}
(\bigcup_{n\in\omega}Y^l_{\alpha,n}(\eta))>0\}.\] 
\quad $\Box$ (Claim 11.1)

\medskip

We now divide the proof of the theorem into two cases.

\medskip

Case 1:\quad $\bar{v}=k$ for a finite $k\in\omega$. 

Fix an $\eta_0\in\; ^{\bar{v}}2$. For any $\alpha\in E$ and $l=0,1$
the fact $\nu_{\bar{v}}(\bigcup_{n\in\omega}C^l_{\alpha,n})=1$ implies that
there exists an $n(\alpha,l)$ such that 
$\eta_0\in C^l_{\alpha,n(\alpha,l)}$.
So for any $m\in\omega$, any distinct $\{\alpha_0,\ldots,\alpha_{m-1}\}
\subseteq E$ and any $h\in\; ^m2$ we have, by Fubini Theorem and
independence,
\[\nu_{\mu}(\bigcap_{i=0}^{m-1}Y^{h(i)}_{\alpha_i,n(\alpha_i,h(i))})\]
\[=\int_{^{\bar{v}}2}\nu_{\mu\smallsetminus\bar{v}}
(\bigcap_{i=0}^{m-1}Y^{h(i)}_{\alpha_i,n(\alpha_i,h(i))}(\eta))d
\nu_{\bar{v}}(\eta)\]
\[=\int_{^{\bar{v}}2}\prod_{i=0}^{m-1}\nu_{\mu\smallsetminus\bar{v}}
(Y^{h(i)}_{\alpha_i,n(\alpha_i,h(i))}(\eta))d\nu_{\bar{v}}(\eta)\]
\[\geqslant (2^{-k})(\prod_{i=0}^{m-1}\nu_{\mu\smallsetminus\bar{v}}
(Y^{h(i)}_{\alpha_i,n(\alpha_i,h(i))}(\eta_0))>0.\]

Case 2:\quad $\bar{v}\geqslant\omega$.\quad
Let \[S\subseteq\lambda\smallsetminus ((\bigcup_{\alpha\in E}v_{\alpha})
\cup\bar{v})\]
be such that $|S|=\bar{v}$. We can factor the forcing $\Bbb P$ to
${\Bbb P}_1*\dot{\Bbb P}_2$ such that 
\[{\cal N}^{\Bbb P}={\cal N}^{{\Bbb P}_1*\dot{\Bbb P}_2},\]
where ${\Bbb P}_1={\cal B}(^{\lambda\smallsetminus S}2)$ and
${\Bbb P}_2=({\cal B}(\;^S2))^{{\cal N}^{{\Bbb P}_1}}$
(see [K2, Theorem 3.13]).
Let $r\in\;^S2$ be a random function  
over ${\cal N}^{{\Bbb P}_1}$. Since $|S|=|\bar{v}|$ we can assume
that $r$ is a random function from $\bar{v}$ to $2$. Since
$C^l_{\alpha,n}\in {\cal N}^{{\Bbb P}_1}$ for any $\alpha\in E$,
$n\in\omega$ and $l=0,1$, and
$\nu_{\bar{v}}(\bigcup_{n\in\omega}C^l_{\alpha,n})=1$,
then there exists an
$n(\alpha,l)$ for any $\alpha\in E$ and $l=0,1$ such that
\[r\in C^l_{\alpha,n(\alpha,l)}.\]
Now for any $m\in\omega$, any distinct $\{\alpha_0,\ldots,\alpha_{m-1}\}
\subseteq E$ and any $h\in\; ^m2$ we have $r\in C$, where
\[C=\bigcap_{i=0}^{m-1}C^{h(i)}_{\alpha_i,n(\alpha_i,h(i))}.\]
This implies $\nu_{\bar{v}}(C)>0$. Hence 
\[\nu_{\mu}(\bigcap_{i=0}^{m-1}Y^{h(i)}_{\alpha_i,
n(\alpha_i,h(i))})\]
\[\geqslant\int_C\nu_{\mu\smallsetminus\bar{v}}
(\bigcap_{i=0}^{m-1}Y^{h(i)}_{\alpha_i,n(\alpha_i,h(i))}(\eta))d
\nu_{\bar{v}}(\eta)\]
\[=\int_C(\prod_{i=0}^{m-1}
\nu_{\mu\smallsetminus\bar{v}}(Y^{h(i)}_{\alpha_i,n(\alpha_i,h(i))}
(\eta)))d\nu_{\bar{v}}(\eta)>0.\]
\quad $\Box$

\medskip

\noindent {\bf Remark:}\quad Theorem 11 is complementary
to Theorem 2 and its corollaries. Combining those results
we conclude that the compactness of a hyperfinite 
Loeb space of size $2^{\omega}$ 
is undecidable under ZFC.

\begin{theorem}

Suppose $\cal M$ is obtained by adding
$\kappa$ random reals to a ZFC model $\cal N$ for some 
regular $\kappa>\omega$. Suppose $\lambda$ is a strong limit cardinal 
such that $cf(\lambda)\leqslant\kappa$. Suppose 
$|\Omega|=\lambda$. Then $\Omega$ is not
compact.

\end{theorem}

\noindent {\bf Proof:}\quad Let $\langle\lambda_{\alpha}:\alpha\in
cf(\lambda)\rangle$ be an increasing sequence such that
$\lambda=\bigcup_{\alpha\in cf(\lambda)}\lambda_{\alpha}$ and
$2^{\lambda_{\alpha}}<\lambda_{\alpha+1}$ for each $\alpha\in cf(\lambda)$.
By similar arguments in the proof of theorem 11 
we can assume that ${\cal B}(\Omega)$ is homogeneous
and $\tau({\cal B}(\Omega))=\lambda$. Note that the cardinality of
any positive measure internal subset of $\Omega$ is $\lambda$.

By Maharam Theorem there is a measure preserving
isomorphism $\Phi$ from ${\cal B}(\Omega)$ to ${\cal B}(^{\lambda}2)$.
Using the same notation as in Theorem 11 
let $A_{\gamma}\subseteq\Omega$ be measurable such that
$\Phi(\bar{A}_{\gamma})=\bar{B}_{\gamma}$ for each $\gamma\in\lambda$. 
By the same argument
as in Claim 10.1 we can find a $cf(\lambda)$-independent sequence
\[\langle C_{\alpha}\subseteq\Omega:
\alpha\in cf(\lambda)\rangle\] of measure $\frac{1}{2}$
such that for any $Z\in [cf(\lambda)]^{cf(\lambda)}$ we have
$\bigcap_{\alpha\in Z}C_{\alpha}=\emptyset$, where $C_{\alpha}=
A_{\gamma_{\alpha}}\Delta A_{\gamma'_{\alpha}}$ for some
different $\gamma_{\alpha},\gamma'_{\alpha}\in 
[\lambda_{\alpha},\lambda_{\alpha+1})$.
Suppose $\Omega$ is compact and assume $\cal C$ is an inner-regular
compact family on $\Omega$. 
Let $X_{\alpha,n}\subseteq C_{\alpha}$ be such that 
\[L_P(\bigcup_{n\in\omega}X_{\alpha,n})=\frac{1}{2}\] and 
let $Y_{\alpha,n}\subseteq\;^{\lambda}2$ 
be Baire sets such that
$\Phi(\bar{X}_{\alpha,n})=\bar{Y}_{\alpha,n}$.
It suffices now to find an $E\in [cf(\lambda)]^{cf(\lambda)}$
and an $n(\alpha)$ for each $\alpha\in E$
such that for any $m\in\omega$ and for any distinct $\{\alpha_0,\ldots,
\alpha_{m-1}\}\subseteq E$ we have 
\[\bigwedge_{i=0}^{m-1}\bar{Y}_{\alpha_i,n(\alpha_i)}\not=0.\]
This is enough to prove the theorem because we have a family
\[\{X_{\alpha,n(\alpha)}:\alpha\in E\}\]
with f.i.p. but
\[\bigcap\{X_{\alpha,n(\alpha)}:\alpha\in E\}\subseteq
\bigcap_{\alpha\in E}C_{\alpha}=\emptyset.\]

For any $\alpha\in cf(\lambda)$ there exists a countable set
$v_{\alpha}\subseteq\kappa$ in $\cal N$ such that
\[\{Y_{\alpha,n}:n\in\omega\}\subseteq 
{\cal N}^{{\cal B}(^{v_{\alpha}}2)}.\]
Choose a $D\in [cf(\lambda)]^{cf(\lambda)}$ such that
\[|\kappa\smallsetminus\bigcup_{\alpha\in D}v_{\alpha}|=\kappa.\]
Let $S\subseteq\kappa\smallsetminus\bigcup_{\alpha\in D}v_{\alpha}$
be such that $|S|=cf(\lambda)$.
Again we factor the forcing ${\Bbb P}={\cal B}(^{\kappa}2)$
to ${\Bbb P}_1*\dot{\Bbb P}_2$ such that 
${\Bbb P}_1={\cal B}(^{\kappa\smallsetminus S}2)$ and
${\Bbb P}_2=({\cal B}(^S2))^{{\cal N}^{{\Bbb P}_1}}$.
Note that $Y_{\alpha,n}\in {\cal N}^{{\Bbb P}_1}$ for any
$\alpha\in D$ and $n\in\omega$. Let
\[R=\bigcup_{\alpha\in D,n\in\omega}supt(Y_{\alpha,n}).\]
Then $R\subseteq\lambda$ and $|R|=cf(\lambda)$.
It is clear that $R$ is unbounded in $\lambda$.
Recall that \[B_{\gamma}=\{f\in\;^{\lambda}2:f(\gamma)=0\}\]
for each $\gamma\in\lambda$ and 
\[\bigvee_{n\in\omega}\bar{Y}_{\alpha,n}=
\bar{B}_{\gamma_{\alpha}}\Delta\bar{B}_{\gamma'_{\alpha}}\]
for each $\alpha\in D$. Notice also that 
$\gamma_{\alpha}, \gamma'_{\alpha}\in R$ for all $\alpha\in D$.
Let $G\subseteq {\Bbb P}_2$ be a ${\cal N}^{{\Bbb P}_1}$ generic
filter. Without loss of generality we assume that
${\Bbb P}_2={\cal B}(^R2)$ since $|S|=|R|$.
Now we define a dense subset $D_{\alpha}\subseteq {\Bbb P}_2$
in ${\cal N}^{{\Bbb P}_1}$ for each $\alpha\in D$. 
For any $Z\subseteq \;^R2$ let 
\[Z^+=\{f\in\;^{\lambda}2:f\!\res\!R\in Z\}\]
and for any $Z\subseteq \;^{\lambda}2$ let 
\[Z^-=\{f\in\;^R2:f=g\!\res\!R\mbox{ for some }g\in Z\}.\]
Define \[D_{\alpha}=\{\bar{Z}:\nu_R(Z)>0\wedge (\exists\beta\in
[\alpha,cf(\lambda))\cap D)(\bar{Z^+}\leqslant\bar{B}_{\gamma_{\beta}}\Delta
\bar{B}_{\gamma'_{\beta}})\}.\]

{\bf Claim 12.1}\quad In ${\cal N}^{{\Bbb P}_1}$, the set 
$D_{\alpha}$ is dense in ${\Bbb P}_2$.

Proof of Claim 12.1:\quad
Given any $\bar{X}\in {\Bbb P}_2$ for some Baire set 
$X\subseteq\;^R2$ with \hfill\break
$\nu_R(X)>0$. Since $supt(X)$ is at most
countable, there exists a $\beta\in [\alpha,cf(\lambda))\cap D$ such that
$\gamma_{\beta},\gamma'_{\beta}\in R\smallsetminus supt(X)$. 
Let $Y=X^+\cap B_{\gamma_{\beta}}\Delta B_{\gamma'_{\beta}}$. Then
\[\nu_R(Y^-)=\nu_R(X)\nu_R((B_{\gamma_{\beta}}\Delta 
B_{\gamma'_{\beta}})^-)>0\] and
$\bar{Y}\leqslant\bar{B}_{\gamma_{\beta}}\Delta\bar{B}_{\gamma'_{\beta}}$.
\quad $\Box$ (Claim 12.1)

\medskip

Let \[E=\{\alpha\in D:
\bar{B}^-_{\gamma_{\alpha}}\Delta\bar{B}^-_{\gamma'_{\alpha}}\in G\}.\]
By Claim 12.1 we have $|E|=cf(\lambda)$.
For each $\alpha\in E$ since 
\[\bigvee_{n\in\omega}\bar{Y}^-_{\alpha,n}=
\bar{B}^-_{\gamma_{\alpha}}\Delta\bar{B}^-_{\gamma'_{\alpha}},\]
then there exists an $n(\alpha)\in\omega$ such that
$\bar{Y}^-_{\alpha,n(\alpha)}\in G$. 
We are done because $G$ is a filter and $supt(Y_{\alpha,n(\alpha)})
\subseteq R$. \quad $\Box$

\medskip

\noindent {\bf Remarks:}\quad (1) Theorem 12 is complementary to
Theorem 6 and Corollary 7.

(2) In this theorem we didn't require 
$\kappa\geqslant 2^{\omega}$ in $\cal N$.

\begin{theorem}

Suppose $\lambda>|V|$, where $V$ is the standard universe,
and $\lambda^{\omega}=\lambda$.
Then there exists a $^*V$, in which every hyperfinite Loeb space $\Omega$
has cardinality $\lambda$ and is not compact.

\end{theorem}

\noindent {\bf Proof:}\quad Note that $|V|=\beth_{\omega}>2^{\omega}$.
Construct a continuous elementary chain of nonstandard universes
\[\langle ^*V_{\alpha}:\alpha\leqslant (2^{\omega})^+\rangle\]
such that for every $\alpha\in (2^{\omega})^+$,

(1) every hyperfinite Loeb space in $^*V_{\alpha}$ has cardinality $\lambda$,
 
(2) $^*V_{\alpha}$ is $\omega_1$-saturated when $\alpha$ is a successor
ordinal.

(3) for any hyperfinite Loeb space $\Omega$ in $^*V_{(2^{\omega})^+}$ if
$\Omega\in\; ^*V_{\alpha}$, then  
$\Omega\cap\; ^*V_{\alpha}$ has Loeb measure zero in $^*V_{\alpha+1}$.

The elementary chain of nonstandard universes satisfying (1), (2) and (3)
exists because at each step one need only to
realize $\leqslant\lambda^{\omega}$ types. We want to show that nonstandard
universe $^*V=\;^*V_{(2^{\omega})^+}$ is the one we want.

Obviously, $^*V_{(2^{\omega})^+}$ is $\omega_1$-saturated.
Suppose there is a compact hyperfinite Loeb space $\Omega$ in $^*V$.
Let $\cal C$ be an inner-regular compact family on $\Omega$.
For every $\alpha\in (2^{\omega})^+$ such that $\Omega\in\; ^*V_{\alpha}$
there is an internal set
$Z_{\alpha}\in\; ^*V_{\alpha+1}$ such that \[L_P(Z_{\alpha})>\frac{1}{2}\]
and \[Z\cap (\Omega\cap\; ^*V_{\alpha})=\emptyset.\]
Now we can find $X_{\alpha,n}\in {\cal C}$ such that
$X_{\alpha,n}\subseteq Z_{\alpha}$ and \[L_P(\bigcup_{n\in\omega}
X_{\alpha,n})=L_P(Z_{\alpha}).\]
Again using the same method as in the proof of Claim 9.1 we can 
find an $E\subseteq (2^{\omega})^+$ with $|E|=(2^{\omega})^+$ 
and an $n(\alpha)\in\omega$ for each $\alpha\in E$ such that
$\{X_{\alpha,n(\alpha)}:\alpha\in E\}$ has f.i.p.
But \[\bigcap_{\alpha\in E}X_{\alpha,n(\alpha)}\subseteq
\bigcap_{\alpha\in E}Z_{\alpha}=\emptyset.\]
\quad $\Box$

\bigskip

We would like to end this section by making a conjecture.

\medskip

\noindent {\bf Conjecture:}\quad It is consistent with 
ZFC that there are no compact hyperfinite Loeb spaces in 
any nonstandard universes.

\medskip

The reader might notice that all results in \S 1 are not
proved by ZFC. But Theorem 8 assumes only ZFC 
plus a consequence of Singular Cardinal Hypothesis. 
So the non-existence of any compact hyperfinite
Loeb spaces would have to violate Singular Cardinal Hypothesis,
which implies the existence of pretty large cardinals.

\bigskip

Department of Mathematics, College of Charleston

Charleston, SC 29424

\medskip

Department of Mathematics, Rutgers University

New Brunswick, NJ 08903

{\em e-mail: rjin@@math.rutgers.edu}

\bigskip

Institute of Mathematics, The Hebrew University

Jerusalem, Israel

\medskip

Department of Mathematics, Rutgers University

New Brunswick, NJ 08903

\medskip

Department of Mathematics, University of Wisconsin

Madison, WI 53706

\bigskip

{\em Sorting:} The first two addresses are the first author's
and the last three are the second author's.

\end{document}